\def\gpic#1{#1
     \medskip\par\noindent{\centerline{\box\graph}} \medskip}
\def\gpic#1{#1
     \medskip\par\noindent{\centerline{\box\graph}} \medskip}
\documentclass[12pt]{article}      
\usepackage{amsmath}
\usepackage{amssymb}
\usepackage{latexsym}
\begin{document}
%  GENERAL MACROS    ++++++++++++++++++++++++++++++++++++++++++++++++++++++

\newcommand{\comment}[1]{}    % type \comment[ put your comment here ]
\newcommand{\hs}{\enspace}
\newcommand{\hhs}{\thinspace}
\newcommand{\real}{\ifmmode {\rm R} \else ${\rm R}$ \fi}
\def\nat{\hbox{\vrule height 7pt width .7pt depth 0pt\hskip -.5pt\bf N}}
	\newcommand{\qed}{\hfill{\setlength{\fboxsep}{0pt}
                  \framebox[7pt]{\rule{0pt}{7pt}}} \newline}

\newcommand{\eqed}{\qquad{\setlength{\fboxsep}{0pt}
                  \framebox[7pt]{\rule{0pt}{7pt}}} }
\newtheorem{theorem}{Theorem}[section]    % omit the percent if you want
\newtheorem{lemma}[theorem]{Lemma}         %   numbering within section.
\newtheorem{corollary}[theorem]{Corollary}
\newtheorem{remark}[theorem]{Remark}
\newtheorem{definition}[theorem]{Definition}
\newtheorem{claim}[theorem]{Claim}
\newtheorem{conjecture}[theorem]{Conjecture}

        % Put all the macros you use in every paper here.
        % I have zillions, so have omitted all but a couple.

\newcommand{\proof }{{\bf Proof: \enspace}}          %  For beginning of proof.
\input epsf.tex

%  PAPER SPECIFIC MACROS   +++++++++++++++++++++++++++++++++++++++++++++++++

% All the macros for this paper only.
\def\gpic#1{#1
%     \midinsert\centerline{\box\graph}\endinsert }
     \medskip\par\noindent{\centerline{\box\graph}} \medskip}
%gpic picture, centered with space

\def\edge{\leftrightarrow}
\def\noedge{\not\leftrightarrow}
\def\twoedge{\Leftrightarrow}
\def\to{\rightarrow}
\def\f{f} 	%to be changed later if required
\def\mn{(K_{1,n-1}, \overline{K_m})}
\def\ca{\lceil\alpha \rceil}
\def\fla{\lfloor\alpha \rfloor}
\def\fra{\{\alpha \}}
\def\n{\notag}
\def\mn{(,K_{1,n-1}, \overline{K_m})}
\def\c{(n-2)^2-2}
\def\dt{\delta}
\def\Dt{\Delta}
\def\ov#1{\overline{#1}}
\def\kn{\ov K_n}
\def\lc{\lceil}
\def\rc{\rceil}
\def\lf{\lfloor}
\def\rf{\rfloor}
\def\ce#1{\par\noindent\centerline{#1}}
\def\skipit#1{}
\def\esub{\subseteq}
\def\sm{\smallskip}

% TITLE INFORMATION  ++++++++++++++++++++++++++++++++++++++++++++++++++++++

\title{Multiple Vertex Coverings by Specified Induced Subgraphs}
\author{ Zolt\'an F\"uredi\thanks{{\it furedi@math-inst.hu} and {\it
z-furedi@math.uiuc.edu}.  Supported in part by the Hungarian National
Science Foundation under grant  OTKA 016389, and by the
National Security Agency under grant MDA904-98-I-0022},
Dhruv Mubayi\thanks{{\it mubayi@math.gatech.edu}. Supported in part by the National Science Foundation
 under grant DMS-9970325.
Current address: School of Mathematics,
Georgia Institute of Technology, Atlanta, GA 30332-0160
},
Douglas B. West\thanks{{\it west@math.uiuc.edu}} \\
University of Illinois \\
Urbana, IL, 61801}
\date{May 25, 1998; revised December, 1999}
\maketitle

% abstract for AMS talk
%
%Given a set $S$ of graphs, we seek the minimum order $f(S)$ of a graph $G$ such
%that each vertex of $G$ belongs to induced subgraphs isomorphic to each graph
%in $S$.  We prove constructively that $f(S)\le 2N-2k$, where $S$ consists of
%$k$ graphs whose orders sum to $N$.  A result of Entringer, Goddard, and
%Henning implies that this bound holds with equality when $S$ consists of a
%clique and an independent set of the same order.  When $H_1$ has minimum degree
%$k$ and $H_2$ is an independent set of order $n$, we prove that
%$f(H_1,H_2)=n+2\sqrt{kn}+O(1)$ as $n\to \infty$.  When $H_1$ is a star of order
%$m$ and $H_2$ is an independent set of order $n$, we compute $f(H_1,H_2)$
%exactly for almost all choices of $m,n$.

 %  ABSTRACT   +++++++++++++++++++++++++++++++++++++++++++++++++++++++

\begin{abstract} Given graphs $H_1, \ldots , H_k$, let $\f(H_1, \ldots, H_k)$ be
the minimum order of a graph $G$ such that for each $i$, the induced copies
of $H_i$ in $G$ cover $V(G)$.  We prove constructively that
$\f(H_1, H_2)\leq 2(n(H_1)+n(H_2)-2)$; equality holds when
$H_1=\ov H_2=K_n$.  We prove that $f(H_1,\kn)=n+2\sqrt{\delta(H_1) n}+O(1)$ as
$n\to \infty$.
We also determine $f(K_{1,m-1},\kn)$ exactly.
%When $H_1=K_{1,m-1}$ and $H_2=\kn$, some exact results are obtained.

\end{abstract}

% MAIN BODY OF PAPER +++++++++++++++++++++++++++++++++++++++++++++++++++
\section{Introduction}

Entringer, Goddard, and Henning \cite{EGH} determined the minimum order of a
simple graph in which every vertex belongs to both a clique of size $m$ and an
independent set of size $n$.  They obtained a surprisingly simple formula for
this value, which they called $f(m,n)$ (an alternative proof using matrix theory
appears in \cite{G1}).

\begin{theorem} \label{EGH} {\rm \cite{EGH}} For $m,n\geq 2$,
$f(m,n)=${\Large $\lceil$}$(\sqrt{m-1}+\sqrt{n-1}~)^2${\Large $\rceil$}.
\end{theorem}

Theorem \ref{EGH} was motivated by a concept introduced
by Chartrand et al.~\cite{CGS} called the {\it framing number}.
A graph $H$ is {\it homogeneously embeddable} in a graph $G$ if, for
all vertices $x\in V(H)$ and $y\in V(G)$, there
exists an embedding of $H$ into $G$ as an induced subgraph
that maps $x$ to $y$.  The
{\it framing number} $fr(H)$ is the minimum order of 
a graph in which $H$ is homogeneously embeddable.
The framing number of a pair of graphs
$H_1$ and $H_2$, written $fr(H_1, H_2)$,
is the minimum order of a graph $G$ in which
both $H_1$ and $H_2$ are homogeneously embeddable.
Thus $fr(K_m,\overline{K_n})=f(m,n)$.  Various
results about the framing number were developed in 
\cite{CGS}.  The framing number of a pair of cycles
is studied in \cite{GHM}. 

When the graphs to be homogeneously embedded are vertex-transitive,
it matters not which vertex of $H$ is mapped to $y\in V(G)$
as long as $y$ belongs to some induced copy of $H$ in $G$.
Determining the framing number for a pair of graphs becomes an
extremal graph covering problem.  We generalize this variation
to more than two graphs.  

\begin{definition} A graph is $(H_1,\ldots,H_k)$-{\em full} if each vertex
belongs to induced subgraphs isomorphic to each of $H_1,\ldots,H_k$.  We use
$\f(H_1,\ldots,H_k)$ to denote the minimum order of an $(H_1,\ldots,H_k)$-full
graph.
\end{definition}

Equivalently, a graph is $(H_1,\ldots,H_k)$-full if for each $i$, the
induced subgraphs isomorphic to $H_i$ cover the vertex set, so we think in
terms of multiple coverings of the vertex set.

Because every vertex in a cartesian product belongs to induced subgraphs
isomorphic to each factor, we have $\f(H_1,\ldots,H_k)\le \prod_i n(H_i)$, where
$n(G)$ denotes the order of $G$.  In fact, $\f(H_1,\ldots,H_k)$ is much smaller.
Our constructions in Section 2 yield $\f(H_1,\ldots,H_k)\le 2\sum_i (n(H_i)-1)$.
Also, if $k-1$ is a prime power and $n(H_i) <k$ for each $i$, 
then  $\f(H_1,\ldots,H_k)\le (k-1)^2$.
By Theorem \ref{EGH}, the first construction is optimal when $k=2$ for $H_1=K_n$
and $H_2=\kn$. We also provide a construction when $H_1$ is arbitrary and
$H_2=\kn$ that is asymptotically sharp up to an additive constant.

In Section 3, we prove a general lower bound in terms of the order of
$H_2$, the maximum degree of $H_2$, and the minimum degree of $H_1$.
In Section 4,  we determine $\f(K_{1,m-1}, \kn)$ exactly (the related
parameter $\f(K_{m,m}, {\ov K_n})$ is studied in \cite{He1}).
In Section 5, we present several open problems.

Since $f(H_1,\ldots,H_k)=f(\ov H_1,\ldots,\ov H_k)$,
all our results yield corresponding results for complementary conditions.
  We note also that there
is an $(H_1,\ldots,H_k)$-full graph for each order exceeding the minimum, since
duplicating a vertex in such a graph yields another $(H_1,\ldots,H_k)$-full
graph.

We consider only simple graphs, denoting the vertex set and edge set of a graph
$G$ by $V(G)$ and $E(G)$, respectively. The {\it order} of $G$
is $n(G)=|V(G)|$.  We use $N_G(v)$ for the neighborhood
of a vertex $v\in V(G)$ (the set of vertices adjacent to $v$), and we let
$N_G[v]=N_G(v) \cup \{v\}$.  The degree of $v$ is $d_G(v)=|N_G(v)|$; we may drop
the subscript $G$.  For $S \subseteq V(G)$, we write $d_S(v)$
for $|N_G(v) \cap S|$.  The {\it independence number} of $G$
is the maximum size of a subset of $V(G)$ consisting of pairwise nonadjacent vertices;
it is denoted by $\alpha(G)$.
  When $S\subseteq V(G)$, we let $N(S)= \bigcup_{v\in S}N(v)$
and let $G[S]$ denote the subgraph induced by $S$. 
% When $v$ is adjacent to all vertices in a set $S$, we write $v\edge S$. 

\section {General Upper Bounds}
Our upper bounds are constructive.

\begin{theorem} \label{realgeneral} 
If $H_1,\ldots,H_k$ are graphs, then
$\f(H_1, \ldots, H_k)\leq 2\sum_{i=1}^k(n(H_i)-1)$.
\end{theorem}
\proof  We construct an $(H_1,\ldots,H_k)$-full graph $G$ with
$2\sum_{i=1}^k(n(H_i)-1)$ vertices.  For $1\le r\le k$, let $H_{r+k}$ be a graph
isomorphic to $H_r$.  For $r\in\{1,\ldots,2k\}$, distinguish a vertex $u_r$ in
$H_r$, and let $N_r=N_{H_r}(u_r)$ and $H_r'=H_r-u_r$.  Construct $G$ from the
disjoint union $H_1'+\dots+H_{2k}'$ by adding, for each $r$, edges making all of
$V(H_r')$ adjacent to all of $N_{r+1}\cup\cdots\cup N_{r+k-1}$, where the
indices are taken modulo $2k$.

By construction, $G$ has the desired order.  For $v\in V(H_r')$ and
$1\le j\le k-1$,
we have $G[v\cup V(H_{r+j}')]\cong H_{r+j}$ (again taking indices modulo $2k$).
Finally, $V(H_r')$ together with any vertex of $V(H_{r-1}')$ induces a copy of
$H_r$ containing $v$.  \qed

Fig. 1 illustrates the construction of Theorem \ref{realgeneral} in the case
$k=2$; an edge to a circle indicates edges to all vertices in the corresponding
set.

\begin{figure}[hbt]
\gpic{
\expandafter\ifx\csname graph\endcsname\relax \csname newbox\endcsname\graph\fi
\expandafter\ifx\csname graphtemp\endcsname\relax \csname newdimen\endcsname\graphtemp\fi
\setbox\graph=\vtop{\vskip 0pt\hbox{%
    \special{pn 8}%
    \special{ar 500 500 217 217 0 6.28319}%
    \special{ar 500 500 500 500 0 6.28319}%
    \graphtemp=.5ex\advance\graphtemp by 0.500in
    \rlap{\kern 0.500in\lower\graphtemp\hbox to 0pt{\hss $N_1$\hss}}%
    \graphtemp=.5ex\advance\graphtemp by 0.756in
    \rlap{\kern 0.756in\lower\graphtemp\hbox to 0pt{\hss  \hss}}%
    \graphtemp=.5ex\advance\graphtemp by 1.080in
    \rlap{\kern 0.386in\lower\graphtemp\hbox to 0pt{\hss $H_1'$\hss}}%
    \special{ar 2000 500 217 217 0 6.28319}%
    \special{ar 2000 500 500 500 0 6.28319}%
    \graphtemp=.5ex\advance\graphtemp by 0.500in
    \rlap{\kern 2.000in\lower\graphtemp\hbox to 0pt{\hss $N_2$\hss}}%
    \graphtemp=.5ex\advance\graphtemp by 0.756in
    \rlap{\kern 1.744in\lower\graphtemp\hbox to 0pt{\hss  \hss}}%
    \graphtemp=.5ex\advance\graphtemp by 0.386in
    \rlap{\kern 1.403in\lower\graphtemp\hbox to 0pt{\hss $H_2'$\hss}}%
    \special{ar 2000 2000 217 217 0 6.28319}%
    \special{ar 2000 2000 500 500 0 6.28319}%
    \graphtemp=.5ex\advance\graphtemp by 2.000in
    \rlap{\kern 2.000in\lower\graphtemp\hbox to 0pt{\hss $N_3$\hss}}%
    \graphtemp=.5ex\advance\graphtemp by 1.744in
    \rlap{\kern 1.744in\lower\graphtemp\hbox to 0pt{\hss  \hss}}%
    \graphtemp=.5ex\advance\graphtemp by 1.420in
    \rlap{\kern 2.114in\lower\graphtemp\hbox to 0pt{\hss $H_1'$\hss}}%
    \special{ar 500 2000 217 217 0 6.28319}%
    \special{ar 500 2000 500 500 0 6.28319}%
    \graphtemp=.5ex\advance\graphtemp by 2.000in
    \rlap{\kern 0.500in\lower\graphtemp\hbox to 0pt{\hss $N_4$\hss}}%
    \graphtemp=.5ex\advance\graphtemp by 1.744in
    \rlap{\kern 0.756in\lower\graphtemp\hbox to 0pt{\hss  \hss}}%
    \graphtemp=.5ex\advance\graphtemp by 2.114in
    \rlap{\kern 1.097in\lower\graphtemp\hbox to 0pt{\hss $H_2'$\hss}}%
    \special{pn 28}%
    \special{pa 1000 500}%
    \special{pa 1783 500}%
    \special{fp}%
    \special{pa 2000 1000}%
    \special{pa 2000 1783}%
    \special{fp}%
    \special{pa 1500 2000}%
    \special{pa 717 2000}%
    \special{fp}%
    \special{pa 500 1500}%
    \special{pa 500 717}%
    \special{fp}%
    \hbox{\vrule depth2.500in width0pt height 0pt}%
    \kern 2.500in
  }%
}%
}
\vskip .5pc
\ce{Fig. 1.  An $(H_1,H_2)$-full graph}
\end{figure}
\bigskip

As mentioned earlier, Theorem \ref{realgeneral} yields sharp upper bounds when
$k=2$ by letting $H_1=K_n$ and $H_2=\kn$. 
In general, as pointed out by a referee, the bounds
can be  off from the optimal by at least a factor of two.  To describe
the construction that improves Theorem \ref{realgeneral} in some cases,
we use resolvable designs.  We phrase the constructions in the language of
hypergraphs.  A hypergraph ${\cal H}=(V, E)$ has vertex set $V$ and edge
set $E$ consisting of subsets of $V$. ${\cal H}$ is $k$-{\it uniform} if
every edge has size $k$, and ${\cal H}$ is $k$-{\it regular} if every vertex
lies in exactly $k$ edges. A {\it matching} $M$ in ${\cal H}$ is a set of
pairwise disjoint edges; $M$ is {\it perfect} if the union of its
elements is $V$.

A {\it Steiner system} $S(n,k, 2)$ is an $n$-vertex $k$-uniform hypergraph
in which every pair of vertices appears together in exactly one edge.
It is {\it resolvable} if the edges can be partitioned into perfect matchings.  
Ray-Chaudhuri and Wilson \cite{RayWil} showed that the trivial necessary
condition $n \equiv k$ (mod $k^2-k$) for the existence of a resolvable
$S(n,k,2)$ is also sufficient when $n$ is sufficiently large compared to $k$. 

\begin{theorem} \label{snk2} If a resolvable Steiner system
$S(n,k-1,2)$ exists and $H_1,\ldots,H_t$ are graphs of order less than $k$,
where $t\le(n-1)/(k-2)$, then $\f(H_1, \ldots, H_t)\leq n$. 
\end{theorem} 
\proof
Duplicating vertices cannot decrease $f$, so we may assume that $n(H_i)=k-1$ for
each $i$.  Let $V$ and $E$ be the vertex set and edge set of the resolvable
Steiner system $S(n,k-1,2)$; we construct a graph $G$ on vertex set $V$.  For
$1\le i\le t$, consider the $i$th perfect matching $M_i$ consisting of edges
$E^i_1, \ldots, E^i_{n/(k-1)}$.  For $j=1, \ldots, n/(k-1)$, add edges within
each $E^i_j$ to make a copy of $H_i$.

Since every pair of vertices lies in only one edge of $S(n,k-1,2)$, this
construction is well defined.  To see that the construction is $H_i$-full,
consider an arbitrary 
$v \in V$. Exactly one of the $t$ edges containing
$v$ belongs to the $i$th matching.  This edge
forms a copy of $H_i$ containing $v$. 
\qed

In the special case when $n=(k-1)^2$, such a resolvable Steiner system
is an {\it affine plane}, denoted ${\cal H}_{k-1}$.  It is well known (see,
\cite[page 672]{HOC} or \cite {VW}, for example)
that an affine plane ${\cal H}_{k-1}$ exists when $k-1$ is a power of a prime.
This yields the following.

\begin{corollary} \label{pp} If  ${\cal H}_{k-1}$
exists and $n(H_i) < k$ for each $i$, then 
$\f(H_1, \ldots, H_k)\leq (k-1)^2$. 
\end{corollary}

When $n(H_i)=k-1$ for each $i$,  Corollary \ref{pp}
improves the bound in Theorem \ref{realgeneral} (asymptotically)
by a factor of two.  When $k=2$ and $H_2=\kn$, a slightly different
construction gives nearly optimal bounds for each $H_1$ as $n\to \infty$. 
In Theorem \ref{generallb}, we shall prove that this construction
is asymptotically optimal.

\begin{theorem} \label{general2}  If $H$ has order $m$ and positive minimum
degree $\delta$, then $f(H,\kn)< n+2\sqrt{\delta n}+2\delta$ when
$n\ge 9\delta (m-\delta-1)^2$. 
\end{theorem}
\proof
Let $x$ be a vertex of minimum degree $\delta$ in $H$.  We construct an
$(H, \kn)$-full graph $G$ in terms of a parameter $r$ that we optimize later.
Let $V(G)=U\cup W$, where $U=U_1\cup\dots\cup U_r$ and $W=W_1\cup\dots\cup W_r$.
Let $W$ be an independent set of size $n-1+s$, where $s=\lc n/(r-1)\rc$.
Let each $W_i$ have size $s-1$ or $s$ (set $|W_r|=s-1$ and put the remaining $n$
vertices equitably into $r-1$ sets).  For each $i$, set $G[U_i] \cong H[N(x)]$,
and make all of $U_i$ adjacent to all of $W_i$.

Each $U_i\cup w$ with $w\in W_i$ induces $N_H[x]$; we add edges to complete
copies of $H$.  Let $m'=m-\delta-1$.  For $j\in\{1,2,3\}$, let $T_j$ consist of
$m'$ vertices, one chosen from each of $U_{(j-1)m'+1},\ldots,U_{jm'}$.  This
requires $r\ge 3m'$.
Add edges within each $T_j$ so that $G[T_j]\cong H-N[x]$.
For each $U_i$ that contains a vertex of $T_j$, add edges from $U_i$ to
$T_{j+1}$ (indices modulo 3 here) so that $G[U_i\cup T_{j+1}]\cong H-x$.
For $3m'+1\le i\le r$, add edges from $U_i$ to $T_1$ so that
$G[U_i\cup T_1]\cong H-x$.  This completes the construction of $G$, as
sketched in Fig.~2; dots represent the vertices of $\bigcup T_j$,
and arrows suggest the edges from $U_i$ to $T_{j+1}$.

\begin{figure}[h]
\gpic{
\expandafter\ifx\csname graph\endcsname\relax \csname newbox\endcsname\graph\fi
\expandafter\ifx\csname graphtemp\endcsname\relax \csname newdimen\endcsname\graphtemp\fi
\setbox\graph=\vtop{\vskip 0pt\hbox{%
    \special{pn 8}%
    \special{pa 1110 3700}%
    \special{pa 2590 3700}%
    \special{pa 2590 3145}%
    \special{pa 1110 3145}%
    \special{pa 1110 3700}%
    \special{fp}%
    \special{pa 1480 3700}%
    \special{pa 1480 3145}%
    \special{fp}%
    \special{pa 1850 3700}%
    \special{pa 1850 3145}%
    \special{fp}%
    \special{pa 2220 3700}%
    \special{pa 2220 3145}%
    \special{fp}%
    \special{pa 1110 2960}%
    \special{pa 2590 2960}%
    \special{pa 2590 2775}%
    \special{pa 1110 2775}%
    \special{pa 1110 2960}%
    \special{fp}%
    \special{pa 1480 2960}%
    \special{pa 1480 2775}%
    \special{fp}%
    \special{pa 1850 2960}%
    \special{pa 1850 2775}%
    \special{fp}%
    \special{pa 2220 2960}%
    \special{pa 2220 2775}%
    \special{fp}%
    \special{pn 28}%
    \special{pa 1295 3145}%
    \special{pa 1295 2960}%
    \special{fp}%
    \special{pa 1665 3145}%
    \special{pa 1665 2960}%
    \special{fp}%
    \special{pa 2035 3145}%
    \special{pa 2035 2960}%
    \special{fp}%
    \special{pa 2405 3145}%
    \special{pa 2405 2960}%
    \special{fp}%
    \graphtemp=.5ex\advance\graphtemp by 2.831in
    \rlap{\kern 1.295in\lower\graphtemp\hbox to 0pt{\hss $\bullet$\hss}}%
    \graphtemp=.5ex\advance\graphtemp by 2.831in
    \rlap{\kern 1.665in\lower\graphtemp\hbox to 0pt{\hss $\bullet$\hss}}%
    \graphtemp=.5ex\advance\graphtemp by 2.831in
    \rlap{\kern 2.035in\lower\graphtemp\hbox to 0pt{\hss $\bullet$\hss}}%
    \graphtemp=.5ex\advance\graphtemp by 2.831in
    \rlap{\kern 2.405in\lower\graphtemp\hbox to 0pt{\hss $\bullet$\hss}}%
    \graphtemp=.5ex\advance\graphtemp by 3.515in
    \rlap{\kern 2.405in\lower\graphtemp\hbox to 0pt{\hss $W_1$\hss}}%
    \graphtemp=.5ex\advance\graphtemp by 3.515in
    \rlap{\kern 1.295in\lower\graphtemp\hbox to 0pt{\hss $W_{m'}$\hss}}%
    \special{pn 8}%
    \special{pa 0 2590}%
    \special{pa 0 1110}%
    \special{pa 555 1110}%
    \special{pa 555 2590}%
    \special{pa 0 2590}%
    \special{fp}%
    \special{pa 0 2220}%
    \special{pa 555 2220}%
    \special{fp}%
    \special{pa 0 1850}%
    \special{pa 555 1850}%
    \special{fp}%
    \special{pa 0 1480}%
    \special{pa 555 1480}%
    \special{fp}%
    \special{pa 740 2590}%
    \special{pa 740 1110}%
    \special{pa 925 1110}%
    \special{pa 925 2590}%
    \special{pa 740 2590}%
    \special{fp}%
    \special{pa 740 2220}%
    \special{pa 925 2220}%
    \special{fp}%
    \special{pa 740 1850}%
    \special{pa 925 1850}%
    \special{fp}%
    \special{pa 740 1480}%
    \special{pa 925 1480}%
    \special{fp}%
    \special{pn 28}%
    \special{pa 555 2405}%
    \special{pa 740 2405}%
    \special{fp}%
    \special{pa 555 2035}%
    \special{pa 740 2035}%
    \special{fp}%
    \special{pa 555 1665}%
    \special{pa 740 1665}%
    \special{fp}%
    \special{pa 555 1295}%
    \special{pa 740 1295}%
    \special{fp}%
    \graphtemp=.5ex\advance\graphtemp by 2.405in
    \rlap{\kern 0.870in\lower\graphtemp\hbox to 0pt{\hss $\bullet$\hss}}%
    \graphtemp=.5ex\advance\graphtemp by 2.035in
    \rlap{\kern 0.870in\lower\graphtemp\hbox to 0pt{\hss $\bullet$\hss}}%
    \graphtemp=.5ex\advance\graphtemp by 1.665in
    \rlap{\kern 0.870in\lower\graphtemp\hbox to 0pt{\hss $\bullet$\hss}}%
    \graphtemp=.5ex\advance\graphtemp by 1.295in
    \rlap{\kern 0.870in\lower\graphtemp\hbox to 0pt{\hss $\bullet$\hss}}%
    \graphtemp=.5ex\advance\graphtemp by 2.405in
    \rlap{\kern 0.259in\lower\graphtemp\hbox to 0pt{\hss $W_{m'+1}$\hss}}%
    \graphtemp=.5ex\advance\graphtemp by 1.295in
    \rlap{\kern 0.259in\lower\graphtemp\hbox to 0pt{\hss $W_{2m'}$\hss}}%
    \special{pn 8}%
    \special{pa 1110 0}%
    \special{pa 2590 0}%
    \special{pa 2590 555}%
    \special{pa 1110 555}%
    \special{pa 1110 0}%
    \special{fp}%
    \special{pa 1480 0}%
    \special{pa 1480 555}%
    \special{fp}%
    \special{pa 1850 0}%
    \special{pa 1850 555}%
    \special{fp}%
    \special{pa 2220 0}%
    \special{pa 2220 555}%
    \special{fp}%
    \special{pa 1110 740}%
    \special{pa 2590 740}%
    \special{pa 2590 925}%
    \special{pa 1110 925}%
    \special{pa 1110 740}%
    \special{fp}%
    \special{pa 1480 740}%
    \special{pa 1480 925}%
    \special{fp}%
    \special{pa 1850 740}%
    \special{pa 1850 925}%
    \special{fp}%
    \special{pa 2220 740}%
    \special{pa 2220 925}%
    \special{fp}%
    \special{pn 28}%
    \special{pa 1295 555}%
    \special{pa 1295 740}%
    \special{fp}%
    \special{pa 1665 555}%
    \special{pa 1665 740}%
    \special{fp}%
    \special{pa 2035 555}%
    \special{pa 2035 740}%
    \special{fp}%
    \special{pa 2405 555}%
    \special{pa 2405 740}%
    \special{fp}%
    \graphtemp=.5ex\advance\graphtemp by 0.869in
    \rlap{\kern 1.295in\lower\graphtemp\hbox to 0pt{\hss $\bullet$\hss}}%
    \graphtemp=.5ex\advance\graphtemp by 0.869in
    \rlap{\kern 1.665in\lower\graphtemp\hbox to 0pt{\hss $\bullet$\hss}}%
    \graphtemp=.5ex\advance\graphtemp by 0.869in
    \rlap{\kern 2.035in\lower\graphtemp\hbox to 0pt{\hss $\bullet$\hss}}%
    \graphtemp=.5ex\advance\graphtemp by 0.869in
    \rlap{\kern 2.405in\lower\graphtemp\hbox to 0pt{\hss $\bullet$\hss}}%
    \graphtemp=.5ex\advance\graphtemp by 0.185in
    \rlap{\kern 1.295in\lower\graphtemp\hbox to 0pt{\hss $W_{2m'+1}$\hss}}%
    \graphtemp=.5ex\advance\graphtemp by 0.185in
    \rlap{\kern 2.405in\lower\graphtemp\hbox to 0pt{\hss $W_{3m'}$\hss}}%
    \special{pn 8}%
    \special{pa 3700 2960}%
    \special{pa 3700 740}%
    \special{pa 3145 740}%
    \special{pa 3145 2960}%
    \special{pa 3700 2960}%
    \special{fp}%
    \special{pa 3700 2590}%
    \special{pa 3145 2590}%
    \special{fp}%
    \special{pa 3700 1110}%
    \special{pa 3145 1110}%
    \special{fp}%
    \special{pn 28}%
    \special{pa 2960 2775}%
    \special{pa 3145 2775}%
    \special{fp}%
    \special{pa 2960 925}%
    \special{pa 3145 925}%
    \special{fp}%
    \special{pn 8}%
    \special{pa 3700 2220}%
    \special{pa 3145 2220}%
    \special{fp}%
    \special{pa 3700 1850}%
    \special{pa 3145 1850}%
    \special{fp}%
    \special{pa 3700 1480}%
    \special{pa 3145 1480}%
    \special{fp}%
    \special{pa 2960 2960}%
    \special{pa 2960 740}%
    \special{pa 2775 740}%
    \special{pa 2775 2960}%
    \special{pa 2960 2960}%
    \special{fp}%
    \special{pa 2960 2590}%
    \special{pa 2775 2590}%
    \special{fp}%
    \special{pa 2960 1110}%
    \special{pa 2775 1110}%
    \special{fp}%
    \special{pa 2960 2220}%
    \special{pa 2775 2220}%
    \special{fp}%
    \special{pa 2960 1850}%
    \special{pa 2775 1850}%
    \special{fp}%
    \special{pa 2960 1480}%
    \special{pa 2775 1480}%
    \special{fp}%
    \special{pn 28}%
    \special{pa 3145 2405}%
    \special{pa 2960 2405}%
    \special{fp}%
    \special{pa 3145 2035}%
    \special{pa 2960 2035}%
    \special{fp}%
    \special{pa 3145 1665}%
    \special{pa 2960 1665}%
    \special{fp}%
    \special{pa 3145 1295}%
    \special{pa 2960 1295}%
    \special{fp}%
    \graphtemp=.5ex\advance\graphtemp by 0.925in
    \rlap{\kern 3.441in\lower\graphtemp\hbox to 0pt{\hss $W_{3m'+1}$\hss}}%
    \graphtemp=.5ex\advance\graphtemp by 2.775in
    \rlap{\kern 3.441in\lower\graphtemp\hbox to 0pt{\hss $W_{r}$\hss}}%
    \special{pn 8}%
    \special{pa 1665 2683}%
    \special{pa 1277 2294}%
    \special{fp}%
    \special{sh 1.000}%
    \special{pa 1316 2359}%
    \special{pa 1277 2294}%
    \special{pa 1342 2333}%
    \special{pa 1316 2359}%
    \special{fp}%
    \special{pa 1277 2294}%
    \special{pa 1018 2035}%
    \special{fp}%
    \special{pa 1018 1665}%
    \special{pa 1406 1277}%
    \special{fp}%
    \special{sh 1.000}%
    \special{pa 1341 1316}%
    \special{pa 1406 1277}%
    \special{pa 1367 1342}%
    \special{pa 1341 1316}%
    \special{fp}%
    \special{pa 1406 1277}%
    \special{pa 1665 1018}%
    \special{fp}%
    \special{pa 1850 1018}%
    \special{pa 1850 2017}%
    \special{fp}%
    \special{sh 1.000}%
    \special{pa 1869 1943}%
    \special{pa 1850 2017}%
    \special{pa 1832 1943}%
    \special{pa 1869 1943}%
    \special{fp}%
    \special{pa 1850 2017}%
    \special{pa 1850 2683}%
    \special{fp}%
    \special{pa 2683 2035}%
    \special{pa 2294 2424}%
    \special{fp}%
    \special{sh 1.000}%
    \special{pa 2359 2384}%
    \special{pa 2294 2424}%
    \special{pa 2333 2358}%
    \special{pa 2359 2384}%
    \special{fp}%
    \special{pa 2294 2424}%
    \special{pa 2035 2683}%
    \special{fp}%
    \hbox{\vrule depth3.700in width0pt height 0pt}%
    \kern 3.700in
  }%
}%
}
\vskip .5pc
\ce{Fig. 2.  Structure of an $(H,\kn)$-full graph}
\end{figure}

To show that $G$ is $(H,\kn)$-full, it suffices to consider $u\in U_i$ and
$w\in W_i$.  By construction, we have $G[\{w\}\cup U_i\cup T_j]\cong H$ for some
$j$.  The vertices of $W-W_i$ together with $u$ or $w$ form an independent set
of size at least $n+s-1-s+1=n$.

We now choose $r$ to minimize the order of $G$, which equals
$n-1+\delta r + \lc n/(r-1) \rc$.  Calculus suggests the choice
$r=\lc\sqrt{n/\delta}~\rc+1$.  This satisfies the requirement that
$r\ge 3m'$ when $n\ge9\dt(m-\delta-1)^2$.  With this value of $r$, the order of
$G$ is at most $n+\delta(2+\sqrt{n/\delta})+\sqrt{\delta n}$, which equals
the bound claimed.  \qed

In the optimized construction, each $|W_i|$ is about $r|U_i|$.
This reflects the use of $W$ to form the large independent set.  When $n$ is
smaller than $9\delta (m')^2$, we still obtain an improvement on Theorem
\ref{realgeneral} by setting $r=3m'$, where $m'=m-\delta-1$.  The resulting
$(H,\kn)$-full graph has order $n-1+\lc n/(3m'-1)\rc+3\delta m'$, which is less
than $2(n+m)$ when $n$ is bigger than about $3\delta m'$.

\section{A Lower Bound}

In this section we prove a lower bound that holds when the
maximum degree of $H_2$ is less than half the minimum degree of $H_1$.

\begin{theorem} \label{generallb} Let $H_1$ and $H_2$ be graphs such that
$H_1$ has minimum degree $\delta$, and $H_2$ has order $n$ and maximum degree
$\Delta$.  If $2\Delta < \delta$, then 
$$f(H_1, H_2)\ge n+ \left\lc2\sqrt{(n+\Dt)(\dt-2\Dt)}\right\rc -(\dt-\Dt).$$
\end{theorem}
\proof
Let $G$ be an $(H_1,H_2)$-full graph, and choose $A\subset V(G)$ such that
$G[A]\cong H_2$.  Let $v$ be a vertex in $V(G)-A$ with the most neighbors
in $A$.  Since $G$ is $(H_1,H_2)$-full, $v$ belongs to a
set $B\subset V(G)$ such that $G[B]\cong H_2$.  Let $C=V(G)-(A\cup B)$.
Let $k=|A-B|$; we obtain a lower bound on $|C|$ in terms of $k$.

Let $e$ be the number of edges with endpoints in both $C$ and $A\cap B$, and let
$d = |N(v)\cap A|$.  Our lower bound on $C$ arises from the computation below.
The first inequality counts $e$ by the $n-k$ endpoints in $A\cap B$; each lies
in a copy of $H_1$ but has at most $2\Delta$ neighbors outside $C$.  The second
inequality counts $e$ by the endpoints in $C$, using the choice of $v$.
For the third inequality, note that $v$ has at most $\Delta$ neighbors in $B$
and then at most $k$ more in $A-B$.
$$(n-k)(\dt-2\Dt)\le e \le d|C|\le(k+\Delta)|C|.$$
Using the resulting lower bound on $|C|$, we have
\begin{align} 
|V(G)|=|A \cup B|+|C| &\ge n+k+\frac{(n-k)(\dt-2\Dt)}{k+\Dt} \n \\
&= n-(\dt-\Dt)+(k+\Dt)+\frac{(n+\Dt)(\dt-2\Dt)}{k+\Dt}. \n
\end{align}
This expression is minimized by $k+\Dt = \sqrt{(n+\Dt)(\dt-2\Dt)}$,
yielding the desired bound.  \qed

\begin{corollary}\label{indlower} If $H_1$ has minimum degree $\dt$,
then $f(H_1, \kn)= n+2\sqrt{\dt n}+O(1)$ as $n\to \infty$.
\end{corollary}
\proof For $\dt>0$, the upper bound follows from Theorem \ref{general2}, while
the lower bound follows by setting $H_2=\kn$ in Theorem \ref{generallb}.
Now suppose that $\dt=0$ and let $m=n(H_1)$.  Let $\alpha(G,v)$ denote the
maximum size of an independent set containing vertex $v$ in a graph $G$.
Let $s=\min_{v\in V(H_1)}\alpha(H_1,v)$.

We claim that $f(H_1,\kn) = n-s+m$ for $n\ge s$.  For the lower bound, let
$u$ be a vertex of $H_1$ such that $s=\alpha(H_1,u)$.  Completing an independent
$n$-set for a vertex playing the role of $u$ in a copy of $H_1$ requires adding
at least $n-s$ vertices to the $m$ vertices of $H_1$. Since
$H_1$ has at least one isolated vertex,  adding these as isolated vertices
yields an $(H_1,\kn)$-full graph, thus proving the upper bound also.
\qed

By taking complements, one immediately obtains the following
corollary.

\begin{corollary} If $\overline{H_1}$ has minimum degree $\dt$,
then $f(H_1, K_n)= n+2\sqrt{\dt n}+O(1)$ as $n\to \infty$.
\end{corollary}

\section{Stars versus Independent Sets}
In this section we determine $\f(H_1, H_2)$ when $H_1$ is a star of order $m$
and $H_2$ is an independent set of order $n$.  Let $S_m=K_{1,m-1}$.
The problem is rather easy when $n<m$.

\begin{claim} \label{smalln}
For $n<m$, $f(S_m,\kn)=n+m-1$, achieved by $K_{n,m-1}$.
\end{claim}
\proof 
The center of an $m$-star must lie in an independent $n$-set
avoiding its neighbors, so $\f(S_m, \kn)\ge n+m-1$ for all $n$.
When $n<m$, the graph $K_{n,m-1}$ is $(S_m,\kn)$-full.  \qed

The problem behaves much differently when $n\ge m$.  First we provide
a construction.

\begin{figure}[h]
\gpic{
\expandafter\ifx\csname graph\endcsname\relax \csname newbox\endcsname\graph\fi
\expandafter\ifx\csname graphtemp\endcsname\relax \csname newdimen\endcsname\graphtemp\fi
\setbox\graph=\vtop{\vskip 0pt\hbox{%
    \special{pn 8}%
    \special{ar 731 1104 224 224 0 6.28319}%
    \special{ar 1552 1104 224 224 0 6.28319}%
    \special{pn 28}%
    \special{pa 955 1104}%
    \special{pa 1328 1104}%
    \special{fp}%
    \graphtemp=.5ex\advance\graphtemp by 1.000in
    \rlap{\kern 0.582in\lower\graphtemp\hbox to 0pt{\hss $\bullet$\hss}}%
    \graphtemp=.5ex\advance\graphtemp by 1.000in
    \rlap{\kern 0.881in\lower\graphtemp\hbox to 0pt{\hss $\bullet$\hss}}%
    \graphtemp=.5ex\advance\graphtemp by 0.358in
    \rlap{\kern 0.358in\lower\graphtemp\hbox to 0pt{\hss $\bullet$\hss}}%
    \graphtemp=.5ex\advance\graphtemp by 0.358in
    \rlap{\kern 0.507in\lower\graphtemp\hbox to 0pt{\hss $\bullet$\hss}}%
    \graphtemp=.5ex\advance\graphtemp by 0.358in
    \rlap{\kern 0.657in\lower\graphtemp\hbox to 0pt{\hss $\bullet$\hss}}%
    \graphtemp=.5ex\advance\graphtemp by 0.358in
    \rlap{\kern 0.806in\lower\graphtemp\hbox to 0pt{\hss $\bullet$\hss}}%
    \graphtemp=.5ex\advance\graphtemp by 0.358in
    \rlap{\kern 0.955in\lower\graphtemp\hbox to 0pt{\hss $\bullet$\hss}}%
    \graphtemp=.5ex\advance\graphtemp by 0.358in
    \rlap{\kern 1.104in\lower\graphtemp\hbox to 0pt{\hss $\bullet$\hss}}%
    \special{pn 8}%
    \special{pa 358 358}%
    \special{pa 582 1000}%
    \special{pa 507 358}%
    \special{pa 582 1000}%
    \special{pa 657 358}%
    \special{fp}%
    \special{pa 806 358}%
    \special{pa 881 1000}%
    \special{pa 955 358}%
    \special{pa 881 1000}%
    \special{pa 1104 358}%
    \special{fp}%
    \graphtemp=.5ex\advance\graphtemp by 1.000in
    \rlap{\kern 1.403in\lower\graphtemp\hbox to 0pt{\hss $\bullet$\hss}}%
    \graphtemp=.5ex\advance\graphtemp by 1.000in
    \rlap{\kern 1.701in\lower\graphtemp\hbox to 0pt{\hss $\bullet$\hss}}%
    \graphtemp=.5ex\advance\graphtemp by 0.358in
    \rlap{\kern 1.254in\lower\graphtemp\hbox to 0pt{\hss $\bullet$\hss}}%
    \graphtemp=.5ex\advance\graphtemp by 0.358in
    \rlap{\kern 1.403in\lower\graphtemp\hbox to 0pt{\hss $\bullet$\hss}}%
    \graphtemp=.5ex\advance\graphtemp by 0.358in
    \rlap{\kern 1.552in\lower\graphtemp\hbox to 0pt{\hss $\bullet$\hss}}%
    \graphtemp=.5ex\advance\graphtemp by 0.358in
    \rlap{\kern 1.701in\lower\graphtemp\hbox to 0pt{\hss $\bullet$\hss}}%
    \graphtemp=.5ex\advance\graphtemp by 0.358in
    \rlap{\kern 1.851in\lower\graphtemp\hbox to 0pt{\hss $\bullet$\hss}}%
    \special{pa 1254 358}%
    \special{pa 1403 1000}%
    \special{pa 1403 358}%
    \special{pa 1403 1000}%
    \special{pa 1552 358}%
    \special{fp}%
    \special{pa 1552 358}%
    \special{pa 1701 1000}%
    \special{pa 1701 358}%
    \special{pa 1701 1000}%
    \special{pa 1851 358}%
    \special{fp}%
    \special{ar 1104 358 896 119 0 6.28319}%
    \graphtemp=.5ex\advance\graphtemp by 1.478in
    \rlap{\kern 0.731in\lower\graphtemp\hbox to 0pt{\hss $\lf r/2\rf$\hss}}%
    \graphtemp=.5ex\advance\graphtemp by 1.478in
    \rlap{\kern 1.552in\lower\graphtemp\hbox to 0pt{\hss $\lc r/2\rc$\hss}}%
    \graphtemp=.5ex\advance\graphtemp by 0.060in
    \rlap{\kern 1.104in\lower\graphtemp\hbox to 0pt{\hss $n-1+k$\hss}}%
    \graphtemp=.5ex\advance\graphtemp by 1.104in
    \rlap{\kern 0.060in\lower\graphtemp\hbox to 0pt{\hss $X$\hss}}%
    \graphtemp=.5ex\advance\graphtemp by 0.358in
    \rlap{\kern 0.060in\lower\graphtemp\hbox to 0pt{\hss $Y$\hss}}%
    \hbox{\vrule depth1.537in width0pt height 0pt}%
    \kern 2.000in
  }%
}%
}
\vskip .5pc
\ce{Fig.~3.  Construction of an $(S_m,\kn)$-full graph.}
\end{figure}

\begin{lemma} \label{constr} 
For $n\ge m\ge 2$,
$$f(S_m,\kn)\le n+\min_k
\max\left\{k+\left\lc\frac{n-1}{k}\right\rc,2m-3-k\right\}.$$
\end{lemma}
\proof
We define a construction $G$ with parameters $r$ and $k$.
Let $V(G)$ be the disjoint union of $X$ and $Y$, where $|X|=r$ and $|Y|=n-1+k$.
Let $G[X]=K_{\lc r/2\rc,\lf r/2\rf}$, and let $Y$ be an independent set.
Give $k$ neighbors in $Y$ to each vertex in $X$, arranged so that $G$ 
is bipartite and has no isolated vertices.

With $k\ge1$, the size chosen for $Y$ ensures that each vertex lies in
an independent $n$-set.  Keeping $G$ bipartite requires $n-1\ge k$.  This
ensures that each vertex of $X$ lies at the center of an induced star of order
$k+1+\lf r/2\rf$.  Thus we require
$$r/2 \ge m-1-k. \qquad\qquad{\rm (A)}$$
Ensuring that the stars cover $Y$ requires
$$(r-1)k\ge n-1.\qquad\qquad~~{\rm (B)}$$

Given $n\ge m\ge 2$, we choose $r,k$ to minimize $n-1+k+r$, the order of $G$.
Rewrite (A) as $r-1 \ge 2m-3-2k$.  Both (A) and (B) impose
lower bounds on $r-1$ in terms of $k,m,n$; we set
$r-1 = \max\{\lc(n-1)/k\rc,2m-3-2k\}$. This yields the one-variable
minimization in the statement of the lemma. \qed

In fact, the construction of Lemma \ref{constr} is optimal for all $n\ge m$.
%We have noted this for $n>1+(4/9)(m-2)^2$, where the construction achieves the
%lower bound from Theorem \ref{EGH}.
We begin the proof of optimality with a lower
bound that differs from the upper bound by at most 1.

\begin{lemma} \label{lower} 
For $n\ge m\ge 2$,
$$f(S_m,\kn)\ge n+\min_d
\max\left\{d-1+\left\lc\frac{n}{d}\right\rc,2m-2-d\right\}.$$
\end{lemma}
\proof
We strengthen the general argument of Theorem~\ref{generallb}.  Let $G$ be an
$(S_m,\kn)$-full graph.  Let $d$ be the maximum of $| N(v)\cap T|$ such that
$v\in V(G)$ and $T$ is an independent $n$-set in $G$.  Let $A$ be an independent
$n$-set and $x$ a vertex such that $|N(x)\cap A|=d$.

As in the proof of Theorem~\ref{generallb}, we choose $B$ to be an independent
$n$-set containing $x$, let $C=V(G)-(A\cup B)$, and let $k$ be the size
of $A-B$.  With $\delta=1$ and $\Delta=0$, the argument applied there to
the edges joining $C$ and $A\cap B$ yields
$$n-k\le d|C|\le k|C|.$$
Since $d\le k$, we obtain $|V(G)|\ge n+d-1+\lc n/d\rc$.

To complete the proof, we must show that $|V(G)|\ge n+2m-2-d$.
As observed in the proof of Claim~\ref{smalln}, $f(S_m,\kn)\ge n+m-1$ always.
Thus we may assume that $d<m-1$.  In proving a lower bound, we may also assume
that $G$ is a minimal $(S_m,\kn)$-full graph.  In particular, if we delete
any edge of $G$, then the resulting graph is not $S_m$-full.
Let $R_1,\dots,R_t$ be a collection of induced stars of order at least $m$
that cover $V(G)$.  By the minimality of $G$, the vertices that are not centers
of these stars form an independent set.  We consider two cases.

\smallskip
{\it Case 1: The centers of $R_1,\ldots,R_t$ form an independent set.}
In this case, $G$ is a bipartite graph with bipartition $X,Y$, where $X$ is
the set of centers of $R_1,\dots,R_t$ and $Y$ is the set of leaves of
$R_1,\dots,R_t$.  By the definition of $d$ and the restriction to $d<m-1$,
we have $|Y|<n$.  Let $x$ be the center of $R_1$, let $I$ be an independent
$n$-set containing $x$, and let $j=|I\cap X|$.  Each vertex of
$I\cap X$ has at least $m-1$ neighbors in $Y-I$.  Since $|Y-I|<n-(n-j)=j$
and there are at least $j(m-1)$ edges from $I\cap X$ to $Y-I$, some $y\in Y-I$
is incident to at least $m-1$ of these edges.  This gives $y$ at least
$m-1>d$ neighbors in $I$, contradicting the choice of $d$.  Thus this case
cannot occur when $d<m-1$.

\smallskip
{\it Case 2. The centers of $R_1,\dots,R_t$ do not form an independent set.}
By the minimality of $G$, each edge of $G$ is needed to complete some induced
star of order at least $m$ centered at one of its endpoints.  We may assume
that the centers $x$ of $R_1$ and $y$ of $R_2$ are adjacent and that $R_1$
needs the edge $xy$ to reach order $m$.  This implies that $y$ is not adjacent
to any leaf of $R_1$.  In particular, the $m-2$ or more additional vertices that
complete $R_2$ are distinct from those in $R_1$, and
$|V(R_1)\cup V(R_2)|\ge 2m-2$.

Now let $I$ be an independent $n$-set containing $x$.  The vertices of
$R_1\cup R_2$ in $I$ are all neighbors of $y$, and hence there are at most
$d$ of them.  Thus $|V(G)|\ge n-d+2m-2$.
\qed

\medskip
When $d$ in the formula of Lemma~\ref{lower} equals $k$ in the formula of
Lemma~\ref{constr}, the resulting values differ by at most one.  A closer look
at the one-variable optimization shows that the lower bound and the upper bound
differ by at most one.

\begin{theorem} \label{optimal} 
For $n\ge m\ge 2$, the construction of Lemma~\ref{constr} is optimal.
\end{theorem}
\proof
We prove that the lower bound of Lemma~\ref{lower} can be improved to
match the upper bound of Lemma~\ref{constr}.

Choose $A,B,C,d,k$ as in the proof of Lemma \ref{lower}.
If $d\le k-1$ or if there are at most $(d-1)|C|$ edges between $C$ and
$A\cap B$, then we obtain $|C|\ge (n-k)/(k-1)$, which yields
$|V(G)|\ge n+k-1+(n-1)/(k-1)$.  Also $2m-2-d\ge 2m-2-k$.
Setting $k'=k-1$ now yields $|V(G)|\ge n+\max\{k'+\lc(n-1)/k'\rc,2m-3-k'\}$.
Hence the construction is optimal unless there is another construction
satisfying $d=k$ and having more than $(d-1)|C|$ edges between $C$ and
$A\cap B$ (thus there is a $z \in C$
with $d_{A \cap B}(z) \ge d$).  More precisely, for every independent set $A$ of size $n$,
every vertex $x \notin A$ with $d_A(x)=d$, and every independent set
$B$ of size $n$ containing $x$, the following holds:
$$B \supseteq A-N(x) \qquad (*)$$

Choose $z \in C$ with  $d_{A \cap B}(z)=d$, and let $B'$ be an
independent set of size $n$ containing $z$.  Letting $(z, A, B')$
play the role of $(x, A, B)$ in  $(*)$ implies that $B' \supseteq  A-N(z)
\supseteq A-B$.  On the other hand, letting  $(z, B, B')$
play the role of $(x, A, B)$ in $(*)$ implies that $B' \supseteq  B-N(z)
\supseteq B-A$. This implies that $(A-B) \cup (B-A)$ is an independent
set, a contradiction.
\qed

\begin{figure}[h]
\gpic{
\expandafter\ifx\csname graph\endcsname\relax \csname newbox\endcsname\graph\fi
\expandafter\ifx\csname graphtemp\endcsname\relax \csname newdimen\endcsname\graphtemp\fi
\setbox\graph=\vtop{\vskip 0pt\hbox{%
    \special{pn 8}%
    \special{pa 818 1364}%
    \special{pa 818 0}%
    \special{pa 2182 0}%
    \special{pa 2182 1364}%
    \special{pa 818 1364}%
    \special{fp}%
    \special{ar 818 682 682 682 1.570796 4.712389}%
    \special{ar 2182 682 682 682 -1.570796 1.570796}%
    \graphtemp=.5ex\advance\graphtemp by 1.091in
    \rlap{\kern 0.477in\lower\graphtemp\hbox to 0pt{\hss $\bullet$\hss}}%
    \graphtemp=.5ex\advance\graphtemp by 0.818in
    \rlap{\kern 0.477in\lower\graphtemp\hbox to 0pt{\hss $\bullet$\hss}}%
    \graphtemp=.5ex\advance\graphtemp by 0.545in
    \rlap{\kern 0.477in\lower\graphtemp\hbox to 0pt{\hss $\bullet$\hss}}%
    \graphtemp=.5ex\advance\graphtemp by 0.273in
    \rlap{\kern 0.477in\lower\graphtemp\hbox to 0pt{\hss $\bullet$\hss}}%
    \graphtemp=.5ex\advance\graphtemp by 0.273in
    \rlap{\kern 2.455in\lower\graphtemp\hbox to 0pt{\hss $\bullet$\hss}}%
    \special{pa 477 1091}%
    \special{pa 2455 273}%
    \special{pa 477 818}%
    \special{pa 2455 273}%
    \special{pa 477 545}%
    \special{pa 2455 273}%
    \special{pa 477 273}%
    \special{fp}%
    \graphtemp=.5ex\advance\graphtemp by 0.682in
    \rlap{\kern 0.477in\lower\graphtemp\hbox to 0pt{\hss $A-B$\hss}}%
    \graphtemp=.5ex\advance\graphtemp by 0.682in
    \rlap{\kern 1.500in\lower\graphtemp\hbox to 0pt{\hss $A\cap B$\hss}}%
    \graphtemp=.5ex\advance\graphtemp by 0.682in
    \rlap{\kern 2.523in\lower\graphtemp\hbox to 0pt{\hss $B-A$\hss}}%
    \graphtemp=.5ex\advance\graphtemp by 0.273in
    \rlap{\kern 2.591in\lower\graphtemp\hbox to 0pt{\hss $x$\hss}}%
    \graphtemp=.5ex\advance\graphtemp by 0.682in
    \rlap{\kern 0.000in\lower\graphtemp\hbox to 0pt{\hss $d$\hss}}%
    \graphtemp=.5ex\advance\graphtemp by 0.682in
    \rlap{\kern 3.000in\lower\graphtemp\hbox to 0pt{\hss $d$\hss}}%
    \special{ar 1500 1705 682 170 0 6.28319}%
    \graphtemp=.5ex\advance\graphtemp by 1.705in
    \rlap{\kern 1.500in\lower\graphtemp\hbox to 0pt{\hss $\bullet$\hss}}%
    \graphtemp=.5ex\advance\graphtemp by 1.023in
    \rlap{\kern 1.193in\lower\graphtemp\hbox to 0pt{\hss $\bullet$\hss}}%
    \graphtemp=.5ex\advance\graphtemp by 1.023in
    \rlap{\kern 1.398in\lower\graphtemp\hbox to 0pt{\hss $\bullet$\hss}}%
    \graphtemp=.5ex\advance\graphtemp by 1.023in
    \rlap{\kern 1.602in\lower\graphtemp\hbox to 0pt{\hss $\bullet$\hss}}%
    \graphtemp=.5ex\advance\graphtemp by 1.023in
    \rlap{\kern 1.807in\lower\graphtemp\hbox to 0pt{\hss $\bullet$\hss}}%
    \special{pa 1193 1023}%
    \special{pa 1500 1705}%
    \special{pa 1398 1023}%
    \special{pa 1500 1705}%
    \special{pa 1602 1023}%
    \special{pa 1500 1705}%
    \special{pa 1807 1023}%
    \special{fp}%
    \graphtemp=.5ex\advance\graphtemp by 1.705in
    \rlap{\kern 1.636in\lower\graphtemp\hbox to 0pt{\hss $z$\hss}}%
    \graphtemp=.5ex\advance\graphtemp by 1.705in
    \rlap{\kern 0.682in\lower\graphtemp\hbox to 0pt{\hss $C$\hss}}%
    \hbox{\vrule depth1.875in width0pt height 0pt}%
    \kern 3.000in
  }%
}%
}
\vskip .5pc
\ce{Fig.~4.  Final proof of the lower bound.}
\end{figure}

It is worth noting what the result of the one-variable optimization is in
terms of $m$ and $n$.  In particular, the construction achieves a lower
bound resulting from Theorem \ref{EGH} when $n>1+(4/9)(m-2)^2$.

\begin{remark} \label{upcor}
If $n>1+(4/9)(m-2)^2$, then $f(S_m,\kn)= n+\lc2\sqrt{n-1}~\rc$.\hfill\break
If $m \le n\le 1+(4/9)(m-2)^2$, then
$f(S_m,\kn)= n+\lc\frac14(3\beta-\sqrt{\beta^2-8})\sqrt{n-1}~\rc$,
where $2m-3=\beta\sqrt{n-1}$ with $\beta>3$.
\end{remark}
\proof
By Theorem \ref{optimal}, it suffices to minimize over $k$ in Lemma
\ref{constr}.  The term $2m-3-k$ is linear.  The term $k+\lc(n-1)/k\rc$ is
minimized when $k=\lc\sqrt{n-1}~\rc$, where it equals $\lc2\sqrt{n-1}~\rc$.
(When $k=\lc\sqrt{n-1}~\rc$, we let $n-1=k^2-r$ with $r<2k-1$; both formulas
yield $2k-1$ when $r\ge k$ and $2k$ when $r<k$.)

When $2m-3-\lc\sqrt{n-1}~\rc \le\lc2\sqrt{n-1}~\rc$, the construction yields
$f(S_m,\kn)\le n+\lc2\sqrt{n-1}~\rc$.
Since every vertex of an induced star belongs to an induced edge, 
Theorem \ref{EGH} yields $f(S_m,\kn)\ge f(K_2,\kn)\ge n+\lc2\sqrt{n-1}~\rc$.

For smaller $n$, the construction is optimized by choosing $x$ so that
$x+(n-1)/x=2m-3-x$ and letting $k=\lf x \rf$.  The number of vertices
is then $2m-3-k$. 
For large $m$ and $n$, we can approximate the result by ignoring integer parts
and defining $\beta$ by $2m-3=\beta\sqrt{n-1}$.  The solution then occurs at
$x=\frac14(\beta+\sqrt{\beta^2-8}~)\sqrt{n-1}$, and we invoke Theorem
\ref{optimal}.
\qed

\section{Open Problems}
We list several open questions.  The first is the most immediately appealing,
suggested by comparing Theorem \ref{EGH} and Theorem \ref{realgeneral}.
\parindent=0in

\sm
1. Among all choices of an $m$-vertex graph $H_1$ and an $n$-vertex graph $H_2$,
is it true that $f(H_1,H_2)$ is maximized when $H_1$ is a clique and $H_2$ is an
independent set?

\sm
2. Let $G$ be an $S_m$-full graph in which the deletion of any edge produces
a graph that is not $S_m$-full.  Is it true that $G$ must be triangle-free?
\footnote{ this has recently been proved positively in \cite{G2}}

%3. Can one prove directly that $f(S_m,\ov K_{n+1})\ge f(S_m,\kn)$?

\sm
3. Among random graphs, what order is needed so that almost every graph is
$(H_1,\dots,H_k)$-full?

\sm
4. Distinguish a root vertex in each of $H_1,\dots,H_k$.
An $(H_1,\dots,H_k)$-{\it root-full} graph is an $(H_1,\dots,H_k)$-full graph
in which each vertex appears as the root in some induced copy of each $H_i$.
Is it possible to bound the minimum order of such a graph (for arbitrary
choice of roots) in terms of 
$f(H_1,\dots,H_k)$?  (suggested by Fred Galvin)

\sm
5. Similarly, one could require induced copies of each $H_i$ so that for each
$v\in V(G)$ and $x\in H_i$, some copy of $H_i$ occurs with $v$ playing the role
of $x$.  The minimum order of such a graph is the framing number
$fr(H_1,\dots,H_k)$.  How large can $fr(H_1,\ldots,H_k)$ be as a function 
of $f(H_1,\dots,H_k)$?  (suggested by Mike Jacobson)

\section{Acknowledgments}

The authors thank the referees for many suggestions, particularly an
idea leading to Theorem \ref{snk2}.

\end{document}